\documentclass[review,5p,times,twocolumn,round,authoryear,square]{elsarticle}

\usepackage{amssymb,amsmath}
\usepackage{graphics,graphicx}
\usepackage{url,theorem}
\usepackage{booktabs}
\usepackage{xcolor}

\usepackage{epstopdf}

\journal{Chemical Engineering Science}

\newcommand{\torol}[1]{}
\newcommand{\irrevreak}[3]{#1{\ }{\overset{#2}{\longrightarrow}}{\ }#3}
\newcommand{\revreak}[4]{#1{\ }{\overset{#2}{\underset{#3}{\rightleftharpoons}}#4}}

\theorembodyfont{\upshape}
\newtheorem{Thm}{Theorem}
\newtheorem{Lem}{Lemma}
\newtheorem{Rem}{Remark}

\newtheorem{Def}{Definition}
\newcommand{\acr}{absolute concentration robustness}
\newcommand{\FHJg}{Feinberg--Horn--Jackson graph}
\newcommand{\ikde}{induced kinetic differential equation}
\newcommand{\ivp}{initial value problem}
\newcommand{\Mma}{\textit{Mathematica}}
\newcommand{\ode}{ordinary differential equation}
\newcommand{\rrc}{reaction rate coefficient}

\newcommand{\AbsoluteRobustness} {\ifmmode {\textbf{\texttt{AbsoluteRobustness}}} \else {\bf\tt AbsoluteRobustness}\fi}
\newcommand{\AtomMatrix}         {\ifmmode {\textbf{\texttt{AtomMatrix}}} \else {\bf\tt AtomMatrix}\fi}
\newcommand{\Bold}               {\ifmmode {\textbf{\texttt{Bold}}} \else {\bf\tt Bold}\fi}
\newcommand{\Break}              {\ifmmode {\textbf{\texttt{Break}}} \else {\bf\tt Break}\fi}
\newcommand{\Compile}            {\ifmmode {\textbf{\texttt{Compile}}} \else {\bf\tt Compile}\fi}
\newcommand{\Continue}           {\ifmmode {\textbf{\texttt{Continue}}} \else {\bf\tt Continue}\fi}
\newcommand{\CUDALink}           {\ifmmode {\textbf{\texttt{CUDALink}}} \else {\bf\tt CUDALink}\fi}
\newcommand{\DirectedEdges}      {\ifmmode {\textbf{\texttt{DirectedEdges}}} \else {\bf\tt DirectedEdges}\fi}
\newcommand{\Do}                 {\ifmmode {\textbf{\texttt{Do}}} \else {\bf\tt Do}\fi}
\newcommand{\EdgeLabeling}       {\ifmmode {\textbf{\texttt{EdgeLabeling}}} \else {\bf\tt EdgeLabeling}\fi}
\newcommand{\ElementaryReactions}{\ifmmode {\textbf{\texttt{ElementaryReactions}}} \else {\bf\tt ElementaryReactions}\fi}
\newcommand{\ExponentialDistribution}{\ifmmode {\textbf{\texttt{ExponentialDistribution}}} \else {\bf\tt ExponentialDistribution}\fi}
\newcommand{\False}              {\ifmmode {\textbf{\texttt{False}}} \else {\bf\tt False}\fi}
\newcommand{\FindFit}              {\ifmmode {\textbf{\texttt{FindFit}}} \else {\bf\tt FindFit}\fi}
\newcommand{\For}                {\ifmmode {\textbf{\texttt{For}}} \else {\bf\tt For}\fi}
\newcommand{\GetProblem}         {\ifmmode {\textbf{\texttt{GetProblem}}} \else {\bf\tt GetProblem}\fi}
\newcommand{\Goto}               {\ifmmode {\textbf{\texttt{Goto}}} \else {\bf\tt Goto}\fi}
\newcommand{\ImageSize}          {\ifmmode {\textbf{\texttt{ImageSize}}} \else {\bf\tt ImageSize}\fi}
\newcommand{\Label}              {\ifmmode {\textbf{\texttt{Label}}} \else {\bf\tt Label}\fi}
\newcommand{\Manipulate}         {\ifmmode {\textbf{\texttt{Manipulate}}} \else {\bf\tt Manipulate}\fi}
\newcommand{\OpenCLLink}         {\ifmmode {\textbf{\texttt{OpenCLLink}}} \else {\bf\tt OpenCLLink}\fi}
\newcommand{\Parallelize}        {\ifmmode {\textbf{\texttt{Parallelize}}} \else {\bf\tt Parallelize}\fi}
\newcommand{\ParallelMap}        {\ifmmode {\textbf{\texttt{ParallelMap}}} \else {\bf\tt ParallelMap}\fi}
\newcommand{\PlotFunction}       {\ifmmode {\textbf{\texttt{PlotFunction}}} \else {\bf\tt PlotFunction}\fi}
\newcommand{\PlotLabel}          {\ifmmode {\textbf{\texttt{PlotLabel}}} \else {\bf\tt PlotLabel}\fi}
\newcommand{\PoissonDistribution}{\ifmmode {\textbf{\texttt{PoissonDistribution}}} \else {\bf\tt PoissonDistribution}\fi}
\newcommand{\RandomChoice}       {\ifmmode {\textbf{\texttt{RandomChoice}}} \else {\bf\tt RandomChoice}\fi}
\newcommand{\RandomVariate}      {\ifmmode {\textbf{\texttt{RandomVariate}}} \else {\bf\tt RandomVariate}\fi}
\newcommand{\ReactionKinetics}   {\ifmmode {\textbf{\texttt{ReactionKinetics}}} \else {\bf\tt ReactionKinetics}\fi}
\newcommand{\RegularExpression}  {\ifmmode {\textbf{\texttt{RegularExpression}}} \else {\bf\tt RegularExpression}\fi}
\newcommand{\Return}             {\ifmmode {\textbf{\texttt{Return}}} \else {\bf\tt Return}\fi}
\newcommand{\ShowVolpertGraph}   {\ifmmode {\textbf{\texttt{ShowVolpertGraph}}} \else {\bf\tt ShowVolpertGraph}\fi}
\newcommand{\Style}              {\ifmmode {\textbf{\texttt{Style}}} \else {\bf\tt Style}\fi}
\newcommand{\Times}              {\ifmmode {\textbf{\texttt{Times}}} \else {\bf\tt Times}\fi}
\newcommand{\True}               {\ifmmode {\textbf{\texttt{True}}} \else {\bf\tt True}\fi}
\newcommand{\VertexLabeling}     {\ifmmode {\textbf{\texttt{VertexLabeling}}} \else {\bf\tt VertexLabeling}\fi}
\newcommand{\While}              {\ifmmode {\textbf{\texttt{While}}} \else {\bf\tt While}\fi}

\newcommand{\R}{\mathbb{R}}
\newcommand{\Rp}{\mathbb{R}^+}

\newcommand{\RM}{\mathbb{R}^M}

\newcommand{\alphab}{\mbox{\boldmath$\alpha$}}
\newcommand{\betab}{\mbox{\boldmath$\beta$}}
\newcommand{\gammab}{\mbox{\boldmath$\gamma$}}

\newcommand{\nulb}{\ifmmode \mathbf{0}\else \textbf{0}\fi}
\newcommand{\oneb}{\ifmmode \mathbf{1}\else \textbf{1}\fi}
\newcommand{\twob}{\ifmmode \mathbf{2}\else \textbf{2}\fi}

\newcommand{\ab}{\ifmmode \mathbf{a}\else \textbf{a}\fi}
\newcommand{\Ab}{\ifmmode \mathbf{A}\else \textbf{A}\fi}
\newcommand{\bb}{\ifmmode \mathbf{b}\else \textbf{b}\fi}
\newcommand{\Bb}{\ifmmode \mathbf{B}\else \textbf{B}\fi}
\newcommand{\cb}{\ifmmode \mathbf{c}\else \textbf{c}\fi}
\newcommand{\Cb}{\ifmmode \mathbf{C}\else \textbf{C}\fi}
\newcommand{\Db}{\ifmmode \mathbf{D}\else \textbf{D}\fi}
\newcommand{\fb}{\ifmmode \mathbf{f}\else \textbf{f}\fi}
\newcommand{\Fb}{\ifmmode \mathbf{F}\else \textbf{F}\fi}
\newcommand{\gb}{\ifmmode \mathbf{g}\else \textbf{g}\fi}
\newcommand{\hb}{\ifmmode \mathbf{h}\else \textbf{h}\fi}
\newcommand{\kb}{\ifmmode \mathbf{k}\else \textbf{k}\fi}
\newcommand{\Kb}{\ifmmode \mathbf{K}\else \textbf{K}\fi}
\newcommand{\Mb}{\ifmmode \mathbf{M}\else \textbf{M}\fi}
\newcommand{\nb}{\ifmmode \mathbf{n}\else \textbf{n}\fi}
\newcommand{\pb}{\ifmmode \mathbf{p}\else \textbf{p}\fi}
\newcommand{\Pb}{\ifmmode \mathbf{P}\else \textbf{P}\fi}
\newcommand{\qb}{\ifmmode \mathbf{q}\else \textbf{q}\fi}
\newcommand{\rb}{\ifmmode \mathbf{r}\else \textbf{r}\fi}
\newcommand{\sbold}{\ifmmode \mathbf{s}\else \textbf{s}\fi}
\newcommand{\vb}{\ifmmode \mathbf{v}\else \textbf{v}\fi}
\newcommand{\wb}{\ifmmode \mathbf{w}\else \textbf{w}\fi}
\newcommand{\xb}{\ifmmode \mathbf{x}\else \textbf{x}\fi}
\newcommand{\Xb}{\ifmmode \mathbf{X}\else \textbf{X}\fi}
\newcommand{\yb}{\ifmmode \mathbf{y}\else \textbf{y}\fi}
\newcommand{\zb}{\ifmmode \mathbf{z}\else \textbf{z}\fi}

\begin{document}
\begin{frontmatter}
\title{{\ReactionKinetics}---A {\it Mathematica} Package with Applications I.\\
Requirements for a Reaction Kinetics Package}
\author[bme,elte]{J. T\'oth\corref{cor1}\fnref{fn1,fn2}}
\ead{jtoth@math.bme.hu}
\author[bme]{A. L. Nagy\fnref{fn1,fn3}}
\ead{nagyal@math.bme.hu}
\author[nw]{D. Papp}
\ead{dpapp@iems.northwestern.edu}

\address[bme]{Department of Analysis,
Budapest University of Technology and Economics,
Egry J. u. 1., Budapest, Hungary, H-1111}

\address[elte]{Laboratory for Chemical Kinetics of the
Institute of Chemistry,
E\"otv\"os Lor\'and University,
Budapest, Hungary}

\address[nw]{Department of Industrial Engineering and Management Sciences,\\
Northwestern University,
2145 Sheridan Road, Room C210, Evanston, IL 60208}

\fntext[fn1]{Partially supported by the Hungarian National Scientific Foundation, No. 84060.}
\fntext[fn2]{This work is connected to the scientific program of the "Development of quality-oriented and harmonized R+D+I strategy and functional model at BME" project. This project is supported by the New Sz\'echenyi Plan (Project ID: T\'AMOP-4.2.1/B-09/1/KMR-2010-0002).}
\fntext[fn3]{Partially supported by the COST Action CM901: Detailed Chemical Kinetic Models for Cleaner Combustion.}

\cortext[cor1]{Corresponding author}
\begin{abstract}
Requirements are formulated for a reaction kinetics package to be useful for an
as wide as possible circle of users and illustrated
with examples using {\ReactionKinetics}, a {\it Mathematica} based package.
\end{abstract}


\end{frontmatter}
\section{Introduction}
In Part I of our paper we formulate the requirements for
a reaction kinetics package to be useful for an
as wide as possible circle of users.
We try to answer the question:
What should an ideal package know?

In Part II \citep{nagypapptoth} we enumerate the major problems arising when writing and using such a package.

Throughout we try to illustrate everything with the present version of the package and the
kinetic examples are mainly taken so as to join to the lectures presented at the workshop MaCKiE 2011
(\url{http://www.mackie2011.uni-hd.de/}).
In many cases the examples will not show all the fine details
necessary to be given when really using our program: a detailed program documentation will be given later.
Furthermore, trivial examples are chosen in some cases to transparently illustrate the functioning of the package.

Also, there is no place to explain the theoretical background, we refer to the literature,
and also to our future, more detailed work.

\section{History}
In the late sixties, early seventies of the last century, at the time when the computer became
accessible for an ordinary scientist, kineticists immediately started to write codes for at least
three different problems: for parsing large sets of reaction steps,
solving induced kinetic differential equations and simulating the stochastic model of chemical reactions,
see an early review: \cite{garfinkelgarfinkelpringgreenchance}.

Nowadays a few existing widely used program with multiple capabilities are e.g.
\begin{description}
\item
[CHEMKIN]
\url{www.sandia.gov/chemkin/index.html},
\item
[XPPAUT]
\url{www.math.pitt.edu/~bard/xpp/xpp.html}
\item
[KINALC]
\url{garfield.chem.elte.hu/Combustion/kinalc.htm}
\end{description}
but there exists an almost infinite number of them (which we are going to review later; see also \cite{tomlinturanyipilling})
with much less capabilities.

\section{Coverage}
Our major governing principle is that the program
should continuously include recent methods developed in most areas of reaction kinetics modeling.
A conference like MaCKiE 2011 is an excellent occasion, a great help to learn newer methods.

We show a few examples from the recent literature.
\subsection{Testing detailed balance}
\subsubsection{Preliminaries and definition}
The long story starts possibly with \cite{wegscheider}, but the first explicit formulation
(without any formula) of the \emph{principle of detailed balance} has been given by
\cite{fowlermilne}: in real thermodynamic equilibrium
all the subprocesses (whatever they mean) should be in dynamic equilibrium separately in such a way that they do not stop but they proceed with the same velocity in both directions.

One can say that this principle means that time is reversible at equilibrium,
that is why the expression \emph{microscopic reversibility} is usually used as a synonym.

The modern formulation of the principle accepted by IUPAC \cite{goldloeningmacnaughtshemi}
essentially means the same: ``The principle of microscopic reversibility at equilibrium states that, in a system at equilibrium, 
any molecular process and the reverse of that process occur, on the average, at the same rate.''

Neither the above document nor the present authors assert that the principle should hold without any further assumptions; 
for us it is an important hypothesis the fulfilment of which should be checked individually in different models.

Now let us consider a reaction consisting only of reversible steps:
\begin{equation}\label{revreac}
\revreak{\sum_{m=1}^M\alpha(m,p)X(m)}{k_p}{k_{-p}}{\sum_{m=1}^M\beta(m,p)X(m)},\quad(p=1,2,\dots,P)
\end{equation}
where we have $P$ pairs of reaction steps, $\alphab,\betab$ are the matrices 
with \emph{molecularities} as elements, their difference is the \emph{stoichiometric matrix}, and
$k_p\quad (p=-P,\dots,-2,-1,1,2,\dots,P)$ are the \rrc s.
Now, the \ikde\ of the above reaction is as follows:
\begin{eqnarray}\label{kinetic}
&&\hspace{-0.75cm}\dot{\cb}(t)=\nonumber\\
&&\hspace{-0.75cm}\sum_{p=1}^P (\betab(\cdot,p)-\alphab(\cdot,p))\left(k_p\cb(t)^{\alphab(\cdot,p)}-k_{-p}\cb(t)^{\betab(\cdot,p)}\right)
\end{eqnarray}
where $\cb(t)\in\RM$ describes the concentrations of the species at time $t$.
\begin{Def}
Assume that $\cb^*\in\left(\Rp\right)^M$ and it makes the right-hand side of \eqref{kinetic} to be equal to zero, furthermore the condition
\begin{equation}\label{dbconds}
k_p\left(\cb^*\right)^{\alphab(\cdot,p)}=k_{-p}\left(\cb^*\right)^{\betab(\cdot,p)}
\end{equation}
is fulfilled for every $p=1,2\ldots,P$. Then we say that reaction \eqref{revreac} is \emph{detailed balanced at the stationary point $\cb^*$}. If the reaction is detailed balanced at all of its positive stationary points, then it is \emph{detailed balanced}.
\end{Def}

\subsubsection{Methods to test detailed balance}
The most na\"{\i}ve way to check conditions \eqref{dbconds} is the direct one: having found one (or all of) the positive stationary points
substitute it (them) into \eqref{dbconds} and see if it holds or not.
Or, it may be enough to solve equations \eqref{dbconds} for positive $\cb^*$'s to get all the candidates 
in which reaction \eqref{revreac} is detailed balanced.
However a nonlinear system of equations is required to be solved, which makes the problem difficult.

Fortunately, there is a more elegant way to look at the problem, namely \cite{feinbergdb} has proved that the following two conditions are necessary and sufficient for the reaction \eqref{revreac} be detailed balanced.\\
\emph{Condition 1}  (circuit conditions)
Suppose that we have chosen an arbitrary spanning forest for the \FHJg\ of reaction \eqref{revreac}.
It is possible to find a set of $P-N+L$ independent circuits, where $N$ denotes the number of vertices (\emph{complexes}):
$$
N:=|\{\alphab(.,r); r=1,2,\dots,R\}\cup\{\betab(.,r); r=1,2,\dots,R\}|,
$$ 
whereas $L$ is the number of connected components (\emph{linkage classes}) of the \FHJg\
(in which all the complexes are written down exactly once and they are connected with reaction arrows). 
For each of these circuits we write an equation which asserts that the product of the reaction rate coefficients in the clockwise direction and counterclockwise direction is equal. Thus we have $P-N+L$ equations for the reaction rate coefficients.\\
Before formulating Condition 2 we need another important definition.
\begin{Def}
The \emph{deficiency}
is the number of complexes $N$
minus the number $L$ of connected components of the \FHJg\
and the number $S$ of independent reaction steps (or, the rank of the stoichiometric matrix).
\end{Def}
\emph{Condition 2}  (spanning forest conditions)
Assume that reaction \eqref{revreac} shows deficiency $\delta$.
Furthermore, assume that the edges of an arbitrarily selected the spanning forest $F$ has been given an orientation.
Then there exists $\delta$ independent non-trivial solutions to the vector equation
\[
\sum_{(i,j)\in F}\left(\betab(\cdot,j)-\alpha(\cdot,i)\right)a(i,j)=0
\]
for the $a(i,j)$ numbers.
With these solutions one can construct the spanning forest conditions which are
\[
\prod_{(i,j)\in F} \tilde{k}_{ij}^{a(i,j)}=\prod_{(i,j)\in F} \tilde{k}_{ji}^{a(i,j)}
\]
where $\tilde{k}(i,j)$'s are the corresponding reaction rate coefficients (associated to the edge $(i,j)$).

\subsubsection{Applications}
Here we present only one example. Take the model of Wegscheider, i.e.
\GetProblem\textbf{\texttt{["Wegscheider"]}}
to obtain
\begin{center}
\textbf{\texttt{$\{$"A" $\leftrightarrow$ "B", 2 "A" $\leftrightarrow$ "A" + "B"$\}$}}.
\end{center}
Now
\begin{flushleft}
\textbf{\texttt{DetailedBalanced[{"Wegscheider"},{k$_1$,k$_{-1}$,k$_2$,k$_{-2}$}]}}\\[0.4em]
results in\\
{\small\color{brown}DetailedBalanced::nocycle: The FHJ graph of the given formal mechanism has no cycle, so we are given only the spanning forest condition(s).}\\[0.4em]
\textbf{\texttt{\{k$_{-1}$*k$_2$ $==$ k$_{-2}$*k$_1$\}}}
\end{flushleft}
providing the condition for the reaction to be detailed balanced.

There are several applications for e.g. in chirality, models of ion channels, combustion theory
in progress.
\subsection{Absolute concentration robustness}
One might look for sufficient conditions to ensure the independence of a stationary concentration coordinate from outer conditions:
does it only depend on reaction rate coefficients and on no initial concentrations? \citep{shinarfeinberg}.

\subsubsection{Motivation}
A simple example follows. Let us consider the reaction
\[
r1 := \left\{\irrevreak{\text{A} + \text{B}}{k_1}{2 \text{B}}, \irrevreak{\text{B}}{k_2}{\text{A}}\right\},
\]
where A may be interpreted as the inactive form of a protein and B as its active form, e. g.
trypsinogen and trypsin respectively.
(Let us remark in passing that this model is an irreversible version of the Wegscheider model above.)
The induced kinetic differential equation of this reaction is given by
\textbf{\texttt{
DeterministicModel[r1, \{k$_1$, k$_2$\}]]
}}
as
\[
a'(t)=-k_1a(t)b(t)+k_2b(t)\quad b'(t)=k_1a(t)b(t)-k_2b(t).
\]
Calculation of the stationary points (under the condition that \(a(0)=a_0, b(0)=b_0\)) proceeds quite naturally:\\
\textbf{\texttt{
StationaryPoints[r1, \{k$_1$, k$_2$\}, \{a$_0$, b$_0$\},\\
Conditions -> \{a$_0$ + b$_0$ > k$_2$/k$_1$\}, \\
Positivity -> True]
}}

The  stationary concentration \(a^*\)of A, being \(\frac{k_2}{k_1}\), is always positive
(together with \(b^*=a_0+b_0-\frac{k_2}{k_1}\), at least if \(a_0+b_0>\frac{k_2}{k_1}\)), and
does not depend on total mass: the reaction shows \emph{absolute concentration robustness} with respect to A.

However, in the reaction
\(
\revreak{A}{k_1}{k_{-1}}{B}
\)
neither A, nor B has this property, both
coordinates of the stationary concentration vector do depend on \(a_0+b_0\) (as one can easily find it without any program):
\(a^*=k_1\frac{a_0+b_0}{k_{-1}+k_1},\quad b^*=k_{-1}\frac{a_0+b_0}{k_{-1}+k_1}.\)

The question is what kind of general conditions can be given to assure absolute concentration robustness.
An easy/to/check set of conditions has been formulated by \cite{shinarfeinberg}.
\subsubsection{Conditions}
\begin{Thm}
Suppose that a reaction endowed with mass action kinetics
has a positive stationary point, and its deficiency is one.
Then, if there exist two complexes lying in nonterminal strong components of the \FHJg\
which only differ in a single species, then the reaction shows \acr\ with respect to this species.
\end{Thm}

The notions in the theorem are widely used and can be found in the paper cited above.
In the case of the simple example above:
\textbf{\texttt{ReactionsData[r1]["deficiency"]}} gives $\tt{\delta=4-2-1=1}.$

We have shown earlier that there exists a positive stationary point.

To see that there exist two nonterminal complexes only differing in a single species
we draw the \FHJg\ of the reaction, and gladly observe that $(\mathrm{A}+\mathrm{B})-\mathrm{A}=\mathrm{B}.$
(To make it clear, complex  $\mathrm{A}+2\mathrm{B}$ and  $\mathrm{B}$ differs in both A and B.)
\begin{center}\begin{figure}
\includegraphics[width=8cm]{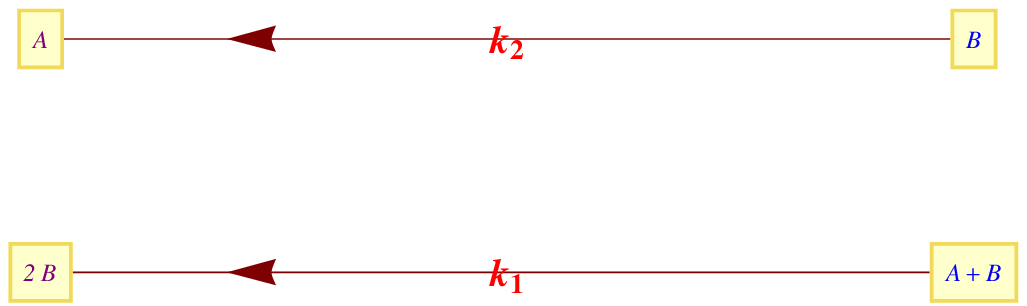}
\caption{The \FHJg\ of the reaction $r1$ obtained by
\textbf{\texttt{
ShowFHJGraph[r1,Style[\#,Red,14]\&/@\{k$_1$, k$_2$\},
VertexLabeling->True,DirectedEdges->True]
}}
}
\end{figure}\end{center}
The last condition is also checked by our program, it would be really hard to test this property in a large reaction set by hand.

Here is how to learn which species of the given reaction are absolutely robust:
\textbf{\texttt{AbsoluteRobustness[r1]}} gives {\tt\{A\}}.

Let us remark that the program is also capable to calculate the stationary concentrations both numerically and \emph{symbolically} in more complicated cases.
Nevertheless one should be cautious.
\begin{Rem}
The theoretical and computational problem lies in showing that a positive stationary point exists.
Why? Because if the \ivp\ to describe a reaction is \(\cb'~=~\fb~\circ~\cb,~\cb(0)~=~\cb_0\)
then the chemically meaningful stationary point \(\cb^*\)
\begin{enumerate}
\item
should obey \(\fb(\cb^*)=\nulb,\)
\item
should be nonnegative,
\item
should be on the level sets of the linear first integrals, and
\item
should be on the level sets of the nonlinear first integrals.
\end{enumerate}
\end{Rem}
The situation is that one cannot in general find (global)
nonlinear first integrals.

Let us consider reaction r1 in detail.
The solutions to \(\fb(\cb^*)=\nulb\) in this case are
\(\left(\begin{array}{c}a^*\\ 0\end{array}\right)\quad (a^*\in\R)\) and
\(\left(\begin{array}{c}\frac{k_2}{k_1}\\ b^*\end{array}\right)\quad (b^*\in\R).\)
These solutions are nonnegative if \(a^*\ge 0\) and \(b^*\ge 0,\) respectively.
They are on the level curves of the linear first integral \(\psi(p,q):=p+q\)
if and only if \(0\le a_0+b_0-\frac{k_2}{k_1}.\)
No nonlinear global first integral---independent from \(\psi\)---exists.

Let us return to detailed balancing.
Certainly, the conditions were not found to treat the trivial examples above.
Their role will be clearer if we cite another---this time more complicated---example from the paper.

Let us consider a model where
ATP is the cofactor in the osmoregulation system of Escherichia coli: EnvZ-OmpR.
\begin{eqnarray*}
r5 = \{\mathrm{X} \leftrightarrow \mathrm{XT} \rightarrow \mathrm{X}_p,
   \mathrm{X}_p + \mathrm{Y} \leftrightarrow \mathrm{X}_p\mathrm{Y} \rightarrow \mathrm{X} + \mathrm{Y}_p,\\
   \mathrm{XT} + \mathrm{Y}_p \leftrightarrow \mathrm{XTY}_p \leftrightarrow \mathrm{XT} + \mathrm{Y}\};
\end{eqnarray*}
The existence of a positive stationary point can be proved either numerically, or
symbolically using \textbf{\texttt{StationaryPoints[r5, Positivity -> True]}} and---waiting for a very long time.

The result given by
\textbf{\texttt{ReactionsData[r5]["deficiency"]}} is
$\tt{\delta=n-l-s=9-3-5=1.}$

\begin{center}
\begin{figure}\label{figshinar}
\includegraphics[width=8cm]{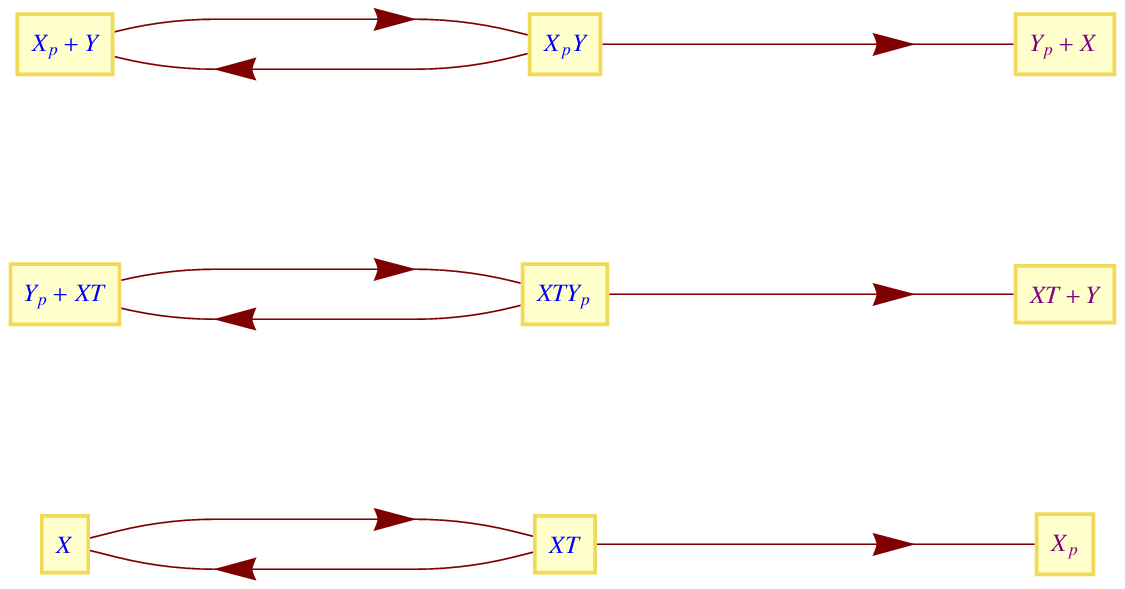}
\caption{The \FHJg\ of reaction $r5$ obtained by\;
\textbf{\texttt{
ShowFHJGraph[r5,VertexLabeling->True,DirectedEdges->True]
}}
with proper complex colourings (purple ones are the terminal complexes).
}
\end{figure}
\end{center}

To see that there exist two nonterminal complexes only differing in a single species
we draw the \FHJg\ of the reaction (Figure \ref{figshinar}) and find that
$(\mathrm{XT} + \mathrm{Y}_p)-\mathrm{XT}=\mathrm{Y}_p,$
therefore the reaction is absolutely robust with respect to $\mathrm{Y}_p,$
the phosphorylated form of the response-regulator,
a species of crucial importance in the reaction.
\subsubsection{The conditions are only sufficient but not necessary}
The fact that the conditions are only sufficient but not necessary can also be seen on the examples given by
\citep{shinarfeinberg}, see also the Supplementary material. Here we give another
example, the one shown by Professor Ross in his lecture (see also \citep[p. 2136]{ross}):
\begin{multline}
jr = \label{reacross}\\
0
\rightarrow \mathrm{X}_1
\leftrightarrow \mathrm{X}_2
\leftrightarrow \mathrm{X}_3
\leftrightarrow \mathrm{X}_4
\leftrightarrow \mathrm{X}_5
\leftrightarrow \mathrm{X}_6
\leftrightarrow \mathrm{X}_7
\leftrightarrow \mathrm{X}_8
\rightarrow 0
\end{multline}

The result given by
\textbf{\texttt{ReactionsData[jr]["deficiency"]}} is
$\delta~=~n~-~l~-~s~=~9~-~1~-~8~=~0,$ (no wonder,
$jr$ is a \emph{compartmental system}),
the theorem cannot be applied,
although all the stationary concentrations are only dependent on the reaction rate coefficients.
Expression
\begin{flushleft}
\textbf{\texttt{
\hspace{-0.2cm}StationaryPoint[jr,\\\{k$_0$,k$_1$,k$_{-1}$,k$_2$,k$_{-2}$,k$_3$,k$_{-3}$,k$_4$,k$_{-4}$,
k$_5$,k$_{-5}$,k$_6$,k$_{-6}$,k$_7$,k$_{-7}$,k$_8$\},\\\{c$_1$[0],c$_2$[0],c$_3$[0],c$_4$[0],c$_5$[0],c$_6$[0],c$_7$[0],c$_8$[0]\}]}}
\end{flushleft}
gives the result symbolically, what we do not reproduce here \emph{verbatim}, because of the length of the formulae.
Here is an equivalent symbolic form:
\[
c_i=
\frac{k_0 \sum _{j=1}^{9-i} \prod _{l=10-j}^8 k_l \left(\prod _{l=i}^{8-j} k_{-l}\right)}{\prod _{l=i}^8 k_l}
\quad (i=1,2,\dots,8)\]
the proof and generalization of which is left to the reader.

\begin{center}
\begin{figure}\label{figross}
\includegraphics[width=9cm]{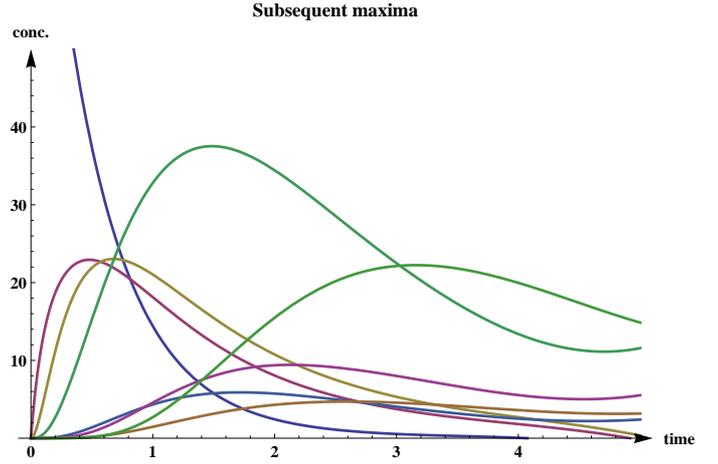}
\caption{Return to the stationary state after a perturbation of the first species
in the reaction (\ref{reacross}) with the \rrc s (\ref{coeffross}).}
\end{figure}
\end{center}

Once here, we can also qualitatively reproduce Fig. 3 of the paper by \cite{ross}.
Let us use the following set of reaction rate coefficients:
\begin{equation}\label{coeffross}
\textbf{\texttt{
rrc = \{0.1,2,0.1,8,5,3,0.4,1,1,6,0.5,4,2,10,1,1\}}}
\end{equation}
and let the initial concentrations be the same as the stationary concentrations except that
we add 100 to the initial value of the initial concentration of X$_1$.
Than the return of the perturbed concentrations
---which was calculated by
\begin{center}
\textbf{\texttt{Concentrations[jr, rrc, ini+100UnitVector[8,1], \{0,4\}]]}}
\end{center}
---
to the (asymptotically stable) original stationary state is as seen on our Fig. \ref{figross}.
More examples and many further interesting details can be found in the mentioned \textit{Science} paper.
\subsection{Improved methods of stochastic simulation}

As in Part II of our paper \cite{nagypapptoth} we have exposed we have built in all the relevant direct and approximate 
methods as well as explicit and implicit ones into our program package. 
For a long while several fancy improvements are known and are involved, here we intend to mention only two ways.

The first one is more theoretical and relies on numerical techniques which have already been applied successfully 
when solving \ode s. 
We mentioned a few of them in \cite{nagypapptoth} but further ideas can and will improve 
the approximation methods of stochastic simulation of chemical reactions.

The other way is based on recent developments and programming tricks. 
The latest versions of \Mma\ enable us to use compiled functions in a more efficient 
and delicate way which are e.g. thousands of times more effective in function evaluation compared to the ''usual'' evaluations 
(see the function \Compile). 
Another direction of recent developments concerns parallel computing (see the functions \Parallelize, \ParallelMap, etc.),
and the use of GPUs (see the functions \CUDALink\ and \OpenCLLink, etc.).

\begin{center}
\begin{figure}\label{figlotka}
\includegraphics[width=9cm]{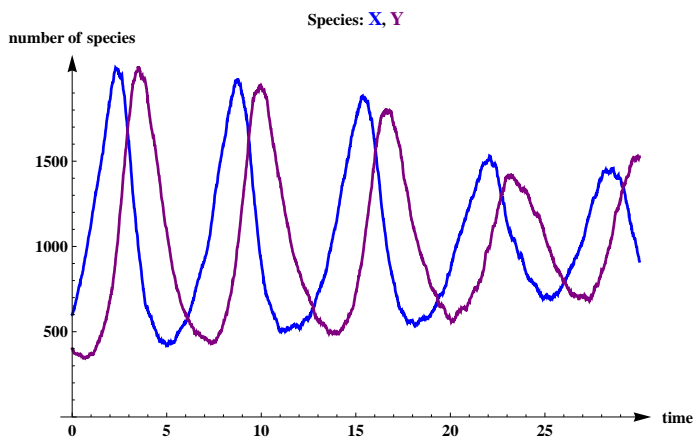}
\caption{Stochastic Lotka--Volterra model with reaction rate coefficients $k_1=1$, $k_2=1/1000$, $k_3=1$ and initial conditions $x_0=600$, $y_0=400$, using tau-leaping methods and compiled functions, where $2.157$ CPU time units were elapsed.}
\end{figure}
\end{center}

Here stands an example: the Lotka--Volterra model is 
\[
\irrevreak{\mathrm{X}}{k_1}{2\mathrm{X}},\,\irrevreak{\mathrm{X}+\mathrm{Y}}{k_2}{2\mathrm{Y}},\,\irrevreak{\mathrm{Y}}{k_3}{0},
\]
and see Figure \ref{figlotka} for the simulation results.

\subsection{Methods of reaction generation and decomposition}
A fundamental problem of stoichiometry is the decomposition of overall reactions into elementary steps.
The definition of ``elementary'' (or ``simple'') reactions varies;
we call here a step elementary if it is of order at most two.
Ideally, a reaction kinetics package should not only provide the user
with the list of possible decompositions
(or at least a large number of decompositions,
as their number may be huge, or even infinite),
but it should also assist the user in the identification of elementary steps themselves,
or even in the generation of possible intermediate species that may take part in an elementary step.
Clearly, the complete automation of this three-step process is a formidable task,
as a large amount of domain specific expertise and data must be incorporated into each step.
A more reasonable goal is to provide the user with a generous list of solutions,
all of which adhere to the basic conservation laws
(e.g., all elementary steps must conserve the number of atoms of each element and total charge)
and other combinatorial constraints imposed by the problem;
and leave all further processing to the user.
At this abstract level the generation of elementary steps with given reactants
and the generation of decompositions of a given overall reaction are essentially equivalent, as we shall see immediately.

Suppose that the species are made of $D$ different elements,
and assign a $(D+1)$-dimensional vector $\ab_m$ to species $X(m)$
where the first $D$ components are the quantities of each constituent,
and the last component is electric charge.
A similar vector can be assigned analogously to every complex as well.
Now, a reaction
\[
\irrevreak{\sum_{m=1}^{M}\alpha_m X(m)}{}{\sum_{m=1}^M\beta_m X(m)}
\]
obeys the laws of atomic and charge balance if and only if the vectors
$\ab_m$ describing the atomic structure of the species satisfy the linear system of equations
\[
\sum_{m=1}^M \alpha_m \ab_{m} = \sum_{m=1}^M\beta_m\ab_{m}.
\]
All combinatorially feasible elementary reactions can be generated by expressing each possible reactant complex (corresponding to a vector 
$\sum_{m=1}^M \alpha_m \ab_{m}$ satisfying $1\leq\sum_m \alpha_{m}\leq 2$) as a linear combination of the vectors $\ab_m$, with nonnegative integer coefficients. Hence, we are looking for the nonnegative integer solutions $\xb$ of a system of linear equations
\[ \Ab\xb = \bb,\]
where $\Ab$ is the \emph{atomic matrix}, with columns $\ab_1$ through $\ab_M$ ($M$ denotes the number of species), and $\bb$ is a vector corresponding to the reactants.
With $M$ species the elementary reactions are obtained by solving $2M+{M \choose 2}$ such systems.

Similarly, if $M$ is the number of species, an $M$-dimensional vector is associated to every reaction involving these species. The $m$th component of the vector shows the change in the quantity of species $X(m)$ in the reaction; if a species is only a reactant, its component is negative, if it is only a product, the coefficient is positive,
if it occurs on both sides, the coefficient might be either positive or negative, or even zero.
The matrix $\gammab$ with column vectors corresponding to the elementary reactions is the stoichiometric matrix of the mechanism. 
A combination of elementary steps with coefficients $\xb$ is a decomposition of the overall reaction if and only if
\[
\gammab\xb = \wb,
\]
where $\wb$ is the vector associated with the overall reaction. Again, we are primarily interested in nonnegative integer solutions $\xb$, though nonnegative rational solutions also can be interpreted as decompositions of the overall reaction.
\torol{
Not surprisingly, this abstract problem has counterparts in perhaps every field concerned with combinatorial structures, from systems biology to computer science and operations research. For instance, \emph{P- and T-invariants of Petri nets}
are solutions of an equivalent problem (albeit with the zero vector on the right-hand side) \citep{martinezsilva}, and so are the \emph{feasible solutions of equality constrained knapsack problems} and the \emph{solutions of linear Diophantine equations} \citep{contejeandevie}.
}

While the complexity of this problem is very high (indeed, deciding whether there exists even a single solution is NP-hard), research in these fields have yielded several relatively practical algorithms and heuristics for its solution.
Further cross-fertilization between fields in this problem is essential, but one must also keep in mind that different applications produce problem instances with strikingly different characteristics, which in turn call for different methods. Indeed, in our investigations in \citep{kovacsvizvaririedeltoth} and \citep{pappvizvari} we have found that entirely different algorithms are the most effective in elementary step generation and in the generation of the decompositions. In particular, note the following obvious differences:
\begin{itemize}
\item in elementary step generation multiple relatively small systems need to be solved, all of which share the coefficient matrix, whereas in decomposition a single large-scale system is solved;
\item in elementary step generation all but one row consist of nonnegative numbers only, in decomposition every row and column may contain entries of different sign;
\item in elementary step generation the number of solutions is always finite, in decomposition it is often infinite.
\end{itemize}

Elaborating on the the last point, if a sequence of steps forms a cycle in which no species are generated or consumed, it can be added to any decomposition an arbitrary number of times. 
A straightforward application of Dickson's lemma \cite{dickson} show is that the converse is also true:
\begin{Lem}
The number of decompositions is finite if and only if the elementary reactions cannot form a cycle. Furthermore, the number of decompositions that do not contain a cycle is always finite.
\end{Lem}
The lemma motivates two further problems related to decompositions: 
first, one needs to be able to decide whether the elementary reactions can form cycles or not, 
which amounts to the solution of a linear programming problem. 
Second, if cycles do exist, the goal changes from generating all decompositions 
to generating all \emph{minimal} (that is, cycle-free) decompositions, 
and all  \emph{minimal cycles} (that is, cycles that cannot be expressed as a sum of two cycles). 
The latter is equivalent to the problem of generating all minimal $P$-invariants of a Petri net
\citep{martinezsilva}.

Another approach to handle the explosion in the number of solutions is to restrict the attention to the simplest decompositions and cycles, those that consist of a small number of steps.

Returning to our discussion of basic requirements for a general purpose reaction kinetics package, it has become clear to us that our package must implement a number of methods for the problems of reaction generation and decomposition, along with an algorithm selection heuristic to choose the one most suited for the problem at hand.

Currently three algorithms are implemented for the solution of these problems, including one of our own developed with very large-scale problems in mind. Additionally, we provide two preprocessing methods to identify the steps that must take part in every decomposition, and the ones that cannot take part in any. (This greatly reduces the computation time for some of the algorithms.) We also provide a heuristic that generates a typically large number of decompositions much faster than the rigorous algorithms, but without any guarantee that it finds all minimal decompositions.

\subsubsection{Example}
We conclude this section with an example for elementary step generation. The oxalate--persulphate--silver oscillator \citep{clarke} involves $16$ species, shown in \mbox{Table \ref{tbl:OPS-species}}; these can form $2\cdot 16+ {16 \choose 2} = 152$ complexes that may be the reactants of an elementary step. The species consist of five atomic constituents, and have charge, hence they can be represented by a $6 \times 16$ atomic matrix, generated by the \AtomMatrix\ command of our package.
From this matrix the \ElementaryReactions\ command generates all elementary steps involving them. In this tiny example each algorithm currently implemented proved to be rather efficient; the number of combinatorially possible elementary steps is $89$.

\begin{table}[bt]
\center
\begin{tabular}{llll}
\toprule
       $\rm Ag^+$ & $\rm Ag^{2+}$ & $\rm H^+$ & $\rm SO_4^{-}$ \\
       $\rm SO_4^{2-}$ & $\rm S_2O_8^{2-}$ & $\rm C_2O_4^{2-}$ & $\rm Ag(C_2O_4)$ \\
       $\rm OH^-$ & $\rm H_2O$ & $\rm CO_2^{-}$ & $\rm O_2$ \\
       $\rm HO_2$ & $\rm H_2O_2$ & $\rm O_2CO_2^-$ & $\rm CO_2$ \\
\bottomrule
\end{tabular}
\caption{Species of the oxalate--persulphate--silver oscillator}
\label{tbl:OPS-species}
\end{table}

In Part II of our paper we shall revisit this problem, and provide a detailed example of obtaining a decomposition for a complex overall reaction.

\section{Applications}
Although there are big differences in the different fields of chemical kinetics, e.g.
the role of thermodynamic data is more important in combustion,
much less important in inorganic chemistry, still we hope the package will be used for more and more realistic applications
in many fields. Here we only mention a few of our previous applications with the ancestors of the present package.
\begin{itemize}
\item Reactions on a surface, reactions in plasma: \cite{sipostotherdi2}.
\item Enzyme kinetics: \cite{tothorlando}.
\item Signal transduction: \cite{tothrospars}.
\item Ion channels: \cite{nagykovacstoth}.
\end{itemize}
\section{Formats}
A useful program should be compatible with usual formats such as CHEMKIN, PrIme, SBML etc.
What we can do at the moment is that we can read in data of different formats but not in an automatic way.
The main method to transform data in different forms is to use pattern matching, including if necessary or useful such tools as \RegularExpression. Suppose we have the file hydrox.dat from the website \url{http://www.math.bme.hu/~nagyal} describing a model of hydrogen oxidation.
Then, we can import and transform it like this:\\
\noindent
{\textbf{\texttt{
Import["hydrox.dat"] /. \{x\_\_\_, {"SPECIES"}, Shortest[y\_\_\_], {"END"}, z\_\_\_\} :> y
}}
The problem is  obviously similar to, or may be considered to be part of parsing.
\section{Speed and Accuracy}
It is important to be able to handle a reaction during a reasonable time interval even if it
consists of a large number of reaction steps and species.
Similarly, the accuracy should also be enough for comparisons with measurements.

Let us mention a few more specific problems where the number of species and of reaction steps can be
really high---beyond the well known areas of atmospheric chemistry, metabolism and combustion.
One might wish to treat reaction chromatography in such a way
that one divides a column into thousands of plates, assumes the same reaction on each of the plates and also assumes (linear)
diffusion between the plates. This is the problem  mentioned by Prof. Trapp in his lecture at the conference.
Let us also mention  here the less known paper by \cite{shapirohorn,shapirohornerr}
which gives a qualitative treatment of such systems using the tools of Chemical Reaction Network Theory.

Or, one might describe the transformation taking place among thousands of polymers
with different molecular weight, as mentioned in the poster by P. van Steenberge \cite{vansteenbergehoogereyniersmarin}.

Another problem might be the treatment of molecules sitting at different energy
levels of which one might have several millions.
A possible first step to treat such a system might be to measure the time needed to solve
the \ikde\ of a model as a function of size as follows.

\noindent
\textbf{\texttt{
lendvay[n\_] :=
 Table[X$_i\leftrightarrow{\ }$X$_{i+1}$, \{i, 1, n - 1\}]\\
Concentrations[lendvay[1000], Array[N[\#] \&, 1998],
     Array[0.1 N[\#] \&, 1000], \{0, 1\}] // Timing
}}

It turned out that such a model can be solved in 30 seconds using \Mma\, without our package,
without \Parallelize\ and without using GPU and similar tools, i. e. one can say, only applying
\emph{Low Performance Computing}.
This is quite a promising start in this direction.

\section{Language: \Mma, what else?}
We do not want to enter into an infinitely long discussion about the advantages and disadvantages of
mathematical program packages, we are only making a few remarks, which
are possibly acceptable by the majority of our readers.

\Mma\ is capable of symbolic
and numeric calculations, creating graphics (in fact, even whole
presentations and publications) within the same framework. In
numerical computations it is not worse
than any other program \citep{weisstein,comparative}.
It uses as many cores/processors as you have in parallel,
it implements OpenCL and CUDA to use GPUs, you can ask it to calculate
the C form of a compiled function if you need, etc.)
Using a symbolic-numeric mathematical programming
language such as \Mma\ is extremely helpful in creating very
transparent programs; code that can even be read as
''pseudo-code'' for those readers who are most reluctant to get
acquainted with the software.

\torol{
On most campuses (in the US, certainly) \Mma\ is available for
students for free through an ''academic headcount'' license to use on
their own computers, and is installed in most computer labs.}

In a meticulous analysis of the consecutive reaction \cite{yablonskyconstalesmartin}
have found a symbolically expressed necessary and sufficient condition for the three
concentration time curves to have a single point of intersection (see also \cite[page 341]{tothsimon}).
As a final application of our package we sow how to solve this 
(and also similar, symbolically untractable problems)
numerically.

\begin{center}\begin{figure}[!h]
\includegraphics[width=8cm]{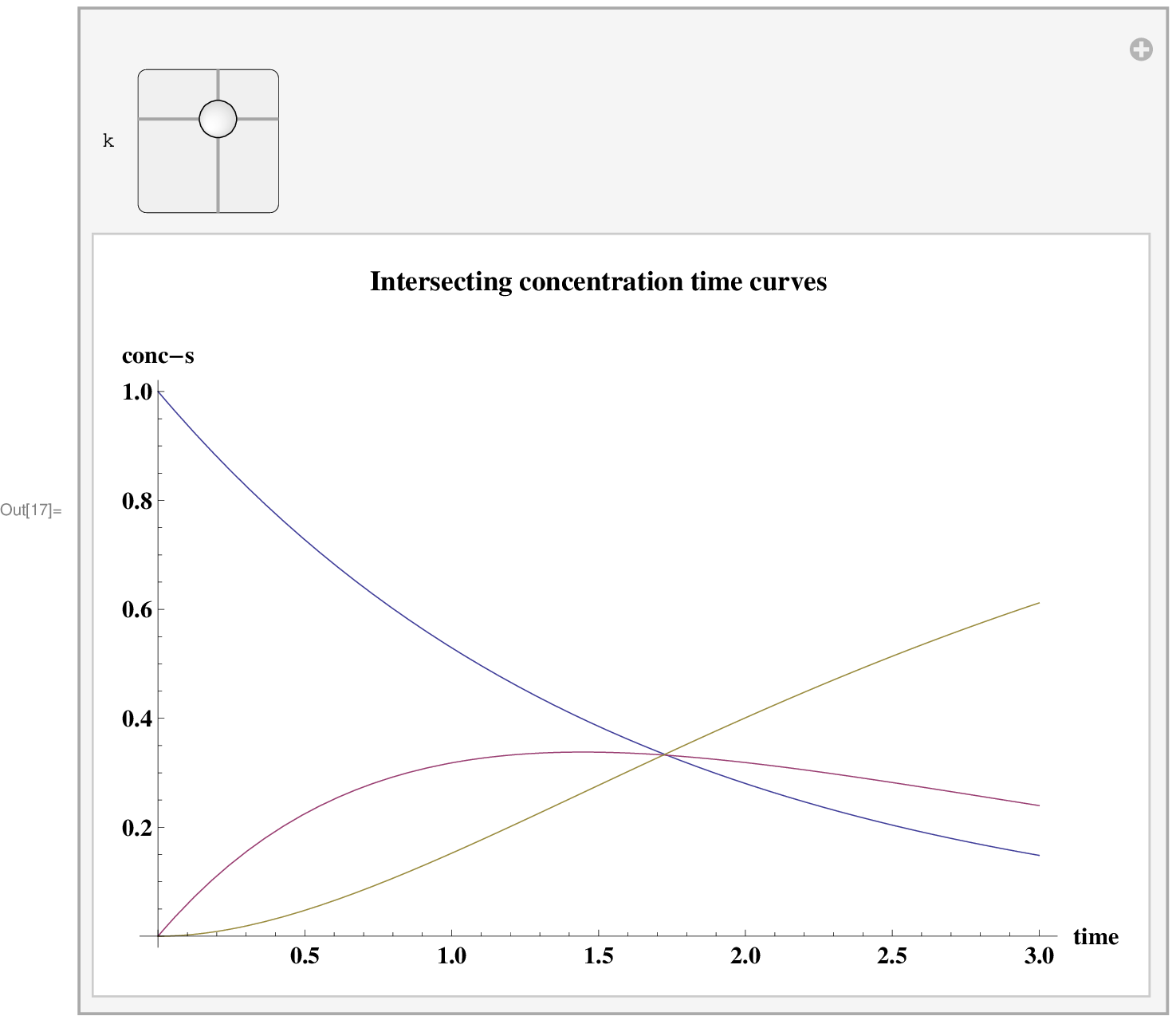}
\caption{Concentration time curves of the consecutive reaction obtained by
\textbf{\texttt{
Manipulate[
 Plot[Evaluate[
   ReplaceAll @@ Concentrations[cons, k, {1, 0, 0}, {0, 3}]], {t, 0, 3}]
}}
}
\end{figure}\end{center}

\section{Outlook}
Having collected so much requirements with illustrations in the second part
we are going to turn to the problems and difficulties when writing and using
a program package pretending to solve so many problems of reaction kinetics.

\bibliography{TothNagyPapp}%
\bibliographystyle{plainnat}%
\end{document}